\input amstex
\documentstyle{amsppt}
\magnification1200
\pagewidth{6.5 true in}
\pageheight{9 true in}
\NoBlackBoxes


\topmatter
\title Sieving and the Erd\H os-Kac theorem
\endtitle
\author Andrew Granville and K. Soundararajan
\endauthor
\abstract We give a relatively easy proof of the Erd\H os-Kac 
theorem via computing moments.  We show how this 
proof extends naturally in a sieve theory context, and 
how it leads to several related results in the 
literature. 
\endabstract
\thanks{Le premier auteur est partiellement soutenu par une bourse
du Conseil  de recherches en sciences naturelles et en g\' enie du
Canada. The second  author is partially supported by the National
Science Foundation.}
\endthanks
\address{D{\'e}partment  de Math{\'e}matiques et Statistique,
Universit{\'e} de Montr{\'e}al, CP 6128 succ Centre-Ville,
Montr{\'e}al, QC  H3C 3J7, Canada}\endaddress
\email{andrew{\@}dms.umontreal.ca}
\endemail
\address{Department of Mathematics, University of Michigan, Ann Arbor,
Michigan 48109, USA} \endaddress \email{ksound{\@}umich.edu}
\endemail
\endtopmatter

\def\phi{\varphi}

\document

\noindent Let $\omega(n)$ denote the number of distinct prime factors 
of the natural number $n$.  The average value of $\omega(n)$ as 
$n$ ranges over the integers below $x$ is 
$$
\frac 1x \sum_{n\leq x} \omega(n) = \frac 1x  \sum_{p\leq x}
 \sum\Sb n\leq x \\ p|n \endSb 1
=   \frac 1x  \sum_{p\leq x}  \left[ \frac xp \right] = \frac 1x  \sum_{p\leq x}  
\left( \frac xp +O(1) \right)  = \log\log x +O(1).
$$
It is natural to ask how $\omega(n)$ is distributed as one varies 
over the integers $n\leq x$.  A famous result of Hardy and Ramanujan [13] tells us 
that $\omega(n)\sim \log \log x$ for almost all $n\leq x$; 
we say that $\omega(n)$ has {\sl normal order} $\log\log n$. 
To avoid confusion let us state this precisely: given $\epsilon>0$ there 
exists $x_\epsilon$ such that if $x\geq x_\epsilon$ is sufficiently large, 
then $(1+\epsilon)\log \log x \geq \omega(n)\geq (1-\epsilon)\log \log x$ 
for all but at most $\epsilon x$ integers $n\leq x$.  
The functions $\log\log n$ and $\log\log x$ are interchangeable here 
since they are very close in value for all but the tiny integers $n\leq x$.

Their proof revolves around the following wonderful inequality which 
they established by induction. Define $\pi_k(x)$
to be the number of integers $n\leq x$ with $\omega(n)=k$.
There exist constants $c_0$, $c_1>0$ such that for any $k\geq 0$ we have
$$
\pi_k(x)<  c_0 \frac{x}{\log x} \ \frac{ (\log\log x \ +c_1)^{k-1}} {(k-1)!} ,
 \tag{1}
$$
for all $x\geq 2$.  Hardy and Ramanujan exploited this by deducing that 
$$
\sum\Sb |k-\log\log x|\geq \epsilon \log\log x\endSb \pi_k(x)
\leq c_0 \frac{x}{\log x} \ \sum\Sb |k-\log\log x|\geq \epsilon \log\log x\endSb
\frac{ (\log\log x \ +c_1)^{k-1}} {(k-1)!} ,
$$
which is easily shown to be about $ x/(\log x)^\alpha$ where $\alpha=\alpha_\epsilon=
\epsilon^2/2+O(\epsilon^3)$, far less than $\epsilon x$. In fact Hardy and Ramanujan 
squeezed a little more out of this idea, showing that if $\kappa(n)\to \infty$ as $n\to \infty$,
no matter how slowly, then
$$
|\omega(n)-\log\log n| \leq \kappa(n) \ \sqrt{\log\log n}
\tag{2}
$$
for almost all integers $n\le x$.

Once we know that $\omega(n)$ has  normal order  $\log\log n$, we can ask finer questions
about the distribution of $\omega(n)$. For instance how is
$\omega(n)-\log\log n$ distributed? More specifically, how big is this typically in absolute value?
Tur\'an [33] found a very simple proof of the Hardy-Ramanujan result by showing that 
$$
\frac 1x \ \sum_{n\leq x} (\omega(n)-\log\log n)^2 = \{ 1+o(1)\} \log\log x .  \tag{3}
$$
One deduces easily that $\omega(n)$ has {\sl normal order} $\log\log n$: 
For, if there are $m_\epsilon(x)$ integers $\leq x$ for which
$|\omega(n)-\log\log n| \geq \epsilon \log\log x$ then by (3), 
$ m_\epsilon(x)\leq (1/\epsilon^2 +o(1)) x/\log\log x $, which is 
$\leq \epsilon x$ for sufficiently large $x$. Indeed the same argument 
also gives (2) for almost all $n\le x$.  

We have now obtained some information about the distribution of $\omega(n)$, its 
average value, and  the average difference between the value and the mean. Next we ask
whether there is a  distribution function for $\omega(n)$? In other words
if, typically, the distance between $\omega(n)$ and $\log\log n$ is roughly of
size $\sqrt{\log\log n}$ can we say anything about the distribution of
$$
\frac{\omega(n)-\log\log n}{\sqrt{\log\log n}}  \ \ ? \tag{4}
$$
In the late 1930s Mark Kac noticed that these developments bore more than a passing
resemblance to developments in probability theory. He suggested that perhaps
this distribution is {\sl normal} and even conjectured certain number theory estimates
which would imply that. Soon after describing this in a lecture, at which Paul Erd\H os 
was in the audience, Erd\H os and Kac were able to announce the result [7]:\ For any $\tau\in \Bbb R$,
the proportion of the integers $n\leq x$ for which $\omega(n)\leq \log\log n + \tau \sqrt{\log\log n}$
tends to the limit
$$
\frac 1{\sqrt{2\pi}} \
\int_{-\infty}^\tau \text{\rm e}^{-t^2/2} dt \tag{5}
$$
as $x\to \infty$. In other words the quantity in (4) is distributed like a 
normal distribution with mean 0 and variance 1.

Erd{\H o}s and Kac's original proof was based on the central limit 
theorem, and Brun's sieve.  A different proof follows from the work 
of Selberg [30] (extending and simplifying the work of [29]) 
who obtained an asymptotic formula for $\pi_k(x)$ 
uniformly in a wide range of $k$.  Yet a third proof is provided 
by Halberstam [11] who showed how to compute the moments 
$$
\sum_{n\le x} (\omega(n)-\log \log x)^k, \tag{6} 
$$
for natural numbers $k$, and showed that these agreed with the 
moments of a normal distribution.  Since the normal distribution is 
well-known to be determined by its moments, he deduced the Erd{\H o}s-Kac 
theorem.  

In this article, we give a simple method to compute the moments 
(6), and in fact we can obtain an asymptotic formula uniformly 
in a wide range of $k$.  Then we discuss how such moments can 
be formulated for more general sequences assuming sieve type 
hypotheses.   

\proclaim{Theorem 1} For any natural number $k$ we 
let $C_k=\Gamma(k+1)/(2^{k/2}\Gamma(k/2+1))$.  
Uniformly for even natural numbers 
$k\le (\log\log x)^{\frac 13}$ we have
$$
\sum_{n\le x} (\omega(n)-\log\log x)^{k} = C_k x(\log\log x)^{k/2} 
\Big(1+O\Big(\frac{k^{\frac 32}}{\sqrt{\log\log x}}\Big)\Big),
$$
and uniformly for odd natural numbers $k \le (\log\log x)^{\frac 13}$ we have
$$
 \sum_{n\le x} (\omega(n)-\log\log x)^{k} \ll C_k x(\log\log x)^{k/2} 
\frac{k^{\frac 32}}{\sqrt{\log\log x}}.
$$
\endproclaim

We will deduce this Theorem from the following technical proposition.  

\proclaim{Proposition 2}  Define 
$$ 
f_p(n) = \cases 1-\frac 1p &\text{if } p|n\\ 
-\frac 1p &\text{if } p\nmid n.\\
\endcases
$$
Let $z \ge 10^6$ be a real number.   Uniformly 
for even natural numbers $k \le (\log\log  z)^{\frac 13}$ we 
have 
$$
\sum_{n\le x} \Big(\sum_{p \le z} f_p(n)\Big)^k 
= C_k x (\log\log z)^{k/2} \Big(1+O\Big(\frac{k^3}{\log\log  z}\Big) \Big) 
+ O(2^k \pi(z)^k), \tag{7a}
$$
while, uniformly for odd natural numbers $k\le (\log\log z)^{\frac{1}{3}}$, 
we have 
$$
\sum_{n\le x} \Big(\sum_{p \le z} f_p(n)\Big)^k \ll C_k x (\log\log z)^{k/2} 
\frac{k^{\frac 32}}{\sqrt{\log\log z} } + 2^{k}\pi(z)^{k}. \tag{7b}
$$
\endproclaim

\demo{Deduction of Theorem 1} 
We seek to evaluate $\sum_{n\le x} (\omega(n)-\log\log x)^{k}$
for natural numbers $k \le (\log\log x)^{\frac 13}$.  
Set $z=x^{\frac 1k}$ and note that, for $n\leq x$,
$$
\omega(n)-\log\log x = \sum_{p\le z} f_p(n) +\sum_{p|n, p>z} 1 + (\sum_{p\le z} 1/p
-\log\log x) = \sum_{p\le z} f_p(n) + O(k).
$$
Thus for some positive constant $c$ we obtain that 
$$
(\omega(n)-\log\log x)^k = \Big(\sum_{p\le z} f_p(n)\Big)^k 
 + O\Big( \sum_{\ell =0}^{k-1} (ck)^{k-\ell} \binom{k}{\ell} 
\Big|\sum_{p\le z} f_p(n)\Big|^{\ell} \Big).
$$
When we sum this up over all integers $n\leq x$ 
the first term above is handled through (7a, b).  To handle 
the remainder terms we need to estimate $\sum_{n\le x} \Big|\sum_{p \le z} 
f_p(n)\Big|^\ell$ for $\ell \le k-1$.  When $\ell$ is even this is 
once again available through (7a).  Suppose $\ell$ is odd.  
By Cauchy-Schwarz we get that 
$$
\sum_{n\le x} \Big|\sum_{p\le z} f_p(n) \Big|^{\ell} 
\le \Big(\sum_{n\le x} \Big(\sum_{p\le z} f_p(n) \Big)^{\ell-1}
\Big)^{\frac 12} 
\Big(\sum_{n\le x} \Big(\sum_{p\le z} f_p(n)\Big)^{\ell+1}\Big)^{\frac 12},
$$
and using (7a) we deduce that this is 
$$
\ll \sqrt{C_{\ell-1}C_{\ell+1}} x (\log\log z)^{\ell/2}.
$$
\enddemo

\demo{Proof of Proposition 2} If $r=\prod_{i} p_i^{\alpha_i}$ is the 
prime factorization of $r$ we put $f_r(n) =\prod_{i} f_{p_i}(n)^{\alpha_i}$. 
Then we may write 
$$
\sum_{n\le x} \Big(\sum_{p\le z} f_p(n)\Big)^k
=\sum_{p_1, \ldots, p_k \le z} \sum_{n\le x} f_{p_1\cdots p_k} (n).
$$
To proceed further, let us consider more generally $\sum_{n\le x} f_r(n)$.  

Suppose $r=\prod_{i=1}^s q_i^{\alpha_i}$ where the $q_i$ are 
distinct primes and $\alpha_i \ge 1$.  Set $R=\prod_{i=1}^{s} q_i$ 
and observe that if $d=(n,R)$ then $f_r(n)=f_r(d)$.  Therefore, with 
$\tau$ denoting the divisor function,  
$$
\align
\sum_{n\le x} f_r(n) &= \sum_{d|R} f_r(d) \sum\Sb n\le x\\ (n,R)=d\endSb 1 
= \sum_{d|R} f_r(d) \Big( \frac{x}{d} \frac{\phi(R/d)}{R/d} + O(\tau(R/d))
\Big) \\
&= \frac{x}{R} \sum_{d|R} f_r(d) \phi(R/d) + O(\tau(R)).
\\
\endalign
$$
Thus seting 
$$
G(r):= \frac 1R \sum_{d|R} f_r(d) \phi(R/d) 
= \prod_{q^\alpha \parallel r} \Big( \frac 1q \Big(1-\frac 1q\Big)^{\alpha} 
+ \Big(\frac{-1}{q}\Big)^{\alpha} \Big(1-\frac 1q\Big)\Big),
$$
we conclude that 
$$
\sum_{n\le x} f_r(n) = G(r) x + O(\tau(R)).
$$

Observe that $G(r)=0$ unless $r$ is square-full and so 
$$
\sum_{n\le x} \Big(\sum_{p\le z}f_p(n)\Big)^k 
= x\sum\Sb p_1, \ldots, p_k \le z\\ p_1\cdots p_k \text{ square-full}\endSb
G(p_1\cdots p_k) + O(2^k \pi(z)^k). \tag{8}
$$
Suppose $q_1 < q_2 <\ldots <q_s$ are the distinct primes in $p_1\cdots p_k$.  
Note that since $p_1\cdots p_k$ is square-full we have $s\le k/2$.  Thus 
our main term above is 
$$
\sum_{s\le k/2} \sum\Sb q_1 < q_2 \ldots < q_{s} \le z \endSb 
\sum\Sb \alpha_1, \ldots, \alpha_s \ge 2\\ \sum_i \alpha_i =k\endSb 
\frac{k!}{\alpha_1! \cdots \alpha_s!} G(q_1^{\alpha_1}\cdots q_s^{\alpha_s}).
$$
When $k$ is even there is a term $s=k/2$ (and all $\alpha_i=2$) 
which gives rise to the Gaussian moments.  This term contributes 
$$
\frac{k!}{2^{k/2}(k/2)!} \sum\Sb q_1,\ldots ,q_{k/2} \le z\\ 
q_i \text{ distinct}\endSb  \ \ \ 
\prod_{i=1}^{k/2} \frac{1}{q_i}\Big(1-\frac{1}{q_i}\Big).  
$$
By ignoring the distinctness condition, we see 
that the sum over $q$'s is bounded above by $(\sum_{p\le z} 
(1-1/p)/p)^{k/2}$.  On the other hand, if we consider $q_1$, $\ldots$, $q_{k/2-1}$ as 
given then the sum over $q_{k/2}$ is plainly at 
least $\sum_{\pi_{k/2} \le p \le z} (1-1/p)/p$ where we 
let $\pi_n$ denote the $n$-th smallest prime.  
Repeating this argument, the sum over the $q$'s is bounded below by 
$(\sum_{\pi_{k/2} \le p \le z} (1-1/p)/p)^{k/2}$.   Therefore the term with $s=k/2$ contributes 
$$
\frac{k!}{(k/2)!2^{k/2}} \Big( \log \log z + O(1+\log \log k) \Big)^{k/2}. 
\tag{9}
$$

To estimate the terms $s<k/2$ we use that $0\le G(q_1^{\alpha_1}\cdots 
q_s^{\alpha_s}) \le 1/(q_1\cdots q_s)$ and so these terms contribute 
$$
\le \sum_{s<k/2}  \frac{k!}{s!} \Big(\sum_{q\le z}\frac 1q\Big)^s 
\sum\Sb \alpha_1,\ldots, \alpha_s \ge 2\\ \sum_i \alpha_i =k \endSb 
\frac{1}{\alpha_1!\cdots \alpha_s!}.
$$
The number of ways of writing $k=\alpha_1+\ldots+\alpha_s$ with 
each $\alpha_i \ge 2$ equals the number of ways of 
writing $k-s=\alpha_1^{\prime}+\ldots+\alpha_s^{\prime}$ where 
each $\alpha_i^{\prime}\ge 1$ and is therefore $\binom{k-s}{s}$.  
Thus these remainder terms contribute 
$$
\le \sum_{s<k/2} \frac{k!}{s!2^s} \binom{k-s}{s} \Big(\log \log z+ O(1)
\Big)^s. \tag{10}
$$
Proposition 2 follows upon combining (8), (9), and (10).  
\enddemo

The main novelty in our proof above is the introduction of the function 
$f_r(n)$ whose expectation over integers $n$ below $x$ is small 
unless $r$ is square-full.  This leads easily to a recognition of the main 
term in the asymptotics of the moments.  Previous approaches expanded 
out $(\omega(n)-\log \log x)^k$ using the binomial theorem, and then there 
are several main terms which must be carefully cancelled out before the 
desired asymptotic emerges.  Our use of this simpler technique was inspired by [25].
Recently Rizwanur Khan [17] builds on this idea to prove that the spacings between 
normal numbers obey a Poisson distribution law.

This technique extends readily to the study of $\omega(n)$ in many other 
sequences.  We formulate this in a sieve like setting:

  Let ${\Cal A} = \{ a_1, \ldots, a_x\}$ be a (multi)-set of $x$ 
(not necessarily distinct) 
natural numbers.  Let ${\Cal A}_d = \#\{ n\le x: \ \ d|a_n\}$.   
We suppose that there is a real valued, non-negative multiplicative function $h(d)$ such that
for square-free $d$ we may write
$$
\Cal A_d = \frac{h(d)}{d} x + r_d. 
$$
It is natural to suppose that $0\le h(d) \le d$ for all square-free $d$, and we do so below.  
Here $r_d$ denotes a remainder term which we expect to be small: either 
small for all $d$, or maybe just small on average over $d$.  

Let $\Cal P$ be any set of primes. In sieve theory one attempts to estimate 
$\# \{ n\leq x:\ (a_n,m)=1\}$ for  $m=\prod_{p\in \Cal P} p$,
in terms of the function $h$ and the error terms $r_d$.
Here we want to understand  the distribution of values of $\omega_{\Cal P} (a)$, 
as  we vary through elements $a$ of $\Cal A$, where 
$\omega_{\Cal P} (a)$ is defined to be the number of primes $p\in \Cal P$ 
which divide $a$.  We  expect that the distribution of $\omega_{\Cal P}(a)$
is normal with  ``mean'' and ``variance'' given by
$$
\mu_{\Cal P}:= \sum_{p\in \Cal P} \frac{h(p)}p\qquad \text{and} \qquad
\sigma_{\Cal P}^2:=  \sum_{p\in \Cal P} \
\frac{h(p)}p\left(1-\frac{h(p)}p\right) \ ,
$$
and wish to find conditions under which this is true. There is a simple
heuristic which explains why this should usually be true:  Suppose that
for each prime $p$ we have a sequence of independent random variables
$b_{1,p}, \dots , b_{x,p}$ each of which is 1 with probability $h(p)/p$
and 0 otherwise; and we let $b_j$ be the product of the primes $p$ for which 
$b_{j,p}=1$. The $b_j$ form a probabilistic model for the $a_j$ satisfying
our sieve hypotheses, the key point being that, in the model, whether or not
$b_j$ is divisible by different primes is independent. One can use the
central limit theorem to show that, as $x\to \infty$, the distribution of
$\omega_{\Cal P}(b)$ becomes normal with mean $\mu_{\Cal P}$ and variance
$\sigma_{\Cal P}^2$.  

\proclaim{Proposition 3}  Uniformly for all
natural numbers $k \le \sigma_{\Cal P}^{\frac 23}$  we have
$$
\sum_{a\in \Cal A} \Big(  \omega_{\Cal P}(a)-\mu_{\Cal
P}\Big)^k = C_k x \sigma_{\Cal P}^k
\Big(1+O\Big(\frac{k^3}{\sigma_{\Cal P}^2}\Big) \Big)
+ O\Big(\mu_{\Cal P}^k  \sum_{d\in D_k(\Cal P)} |r_d|\Big),
$$
if $k$ is even, and 
$$
 \sum_{a\in \Cal A} \Big(  \omega_{\Cal
P}(a)-\mu_{\Cal P}\Big)^k \ll C_k x \sigma_{\Cal P}^k\
\frac{k^{\frac 32}}{\sigma_{\Cal P}} + \mu_{\Cal P}^k  \sum_{d\in
D_k(\Cal P)} |r_d|, 
$$
if $k$ is odd.  Here $D_k(\Cal P)$ denotes the set of squarefree integers which are
the product of at most $k$ primes all from the set $\Cal P$.
\endproclaim
\demo{Proof}  The proof is similar to that of Proposition 2, and so we record 
only the main points.  We define $f_p(a)=1-h(p)/p$ if $p|a$ and $-h(p)/p$ if 
$p\nmid a$.  If $r=\prod_{i} p_i^{\alpha_i}$ is the prime factorization of $r$ 
we put $f_r(a)= \prod_{i} f_{p_i}(a)^{\alpha_i}$.  
Note that $\omega_{\Cal P}(a)-\mu_{\Cal P} = \sum_{p\in {\Cal P}} f_p(a)$, 
and so 
$$
\sum_{a\in {\Cal A}} \Big( \omega_{\Cal P}(a) - \mu_{\Cal P} \Big)^k 
= \sum\Sb p_1, \ldots, p_k \in {\Cal P}\endSb 
\sum_{a\in {\Cal A}} f_{p_1\cdots p_k} (a). \tag{11}
$$

As in Proposition 2, consider more generally $\sum_{a \in {\Cal A}} f_r(a)$.  
Suppose $r=\prod_{i=1}^s q_i^{\alpha_i}$ where the $q_i$ are
distinct primes and each $\alpha_i \ge 1$.  Set $R=\prod_{i=1}^{s}
q_i$ and observe that if $d=(a,R)$ then $f_r(a)=f_r(d)$.  Note that
$$
\align
\sum\Sb a\in \Cal A\\ (a,R)=d\endSb 1
&= \sum_{a\in \Cal A} \  \sum\Sb e|(R/d)\\ de|n\endSb  \mu(e) =
\sum_{e|R/d} \mu(e) \Cal A_{de} \\
&=
x \ \frac {h(d)}d \ \prod_{p|(R/d)} \Big(1-\frac {h(p)}p\Big) + \sum_{e|(R/d)} \mu(e) r_{de} .\\
\endalign
$$
Therefore
$$
\align
\sum_{a\in \Cal A} f_r(a) &= \sum_{d|R} f_r(d) \sum\Sb a\in \Cal A\\ (a,R)=d\endSb 1\\
&= x \sum_{d|R} f_r(d)\ \frac {h(d)}d \  \prod_{p|(R/d)} \Big(1-\frac {h(p)}p\Big) +
\sum_{d|R} f_r(d) \sum_{e|(R/d)} \mu(e) r_{de}
 \\
&=  G(r)x + \sum_{m|R} r_m E(r,m), \tag{12a}\\
\endalign
$$
where
$$
G(r)= \prod\Sb q^\alpha \parallel r \endSb  \Big( \frac
{h(q)}q \Big(1-\frac {h(q)}q\Big)^{\alpha}
+ \Big(\frac{-{h(q)}}{q}\Big)^{\alpha} \Big(1-\frac {h(q)}q\Big)\Big),\tag{12b}
$$
and
$$
E(r,m)= \prod\Sb q^\alpha \parallel r,\ q|m  \endSb  \Big(
\Big(1-\frac {h(q)}q\Big)^{\alpha} -
\Big(\frac{-{h(q)}}{q}\Big)^{\alpha} \Big) \prod\Sb q^\alpha
\parallel r, \ q |(R/m)  \endSb 
\Big(\frac{-{h(q)}}{q}\Big)^{\alpha} . \tag{12c}
$$

We input the above analysis in (11).  Consider first the main terms that arise.  
Notice that $G(r)=0$ unless $r$ is square-full, and so the main terms 
look exactly like the corresponding main terms in Proposition 2.  We record the 
only small difference from the analysis there.  When $k$ is even there is a leading 
contribution from the terms with $s=k/2$ and all $\alpha_i=2$ (in notation analogous 
to Proposition 2); this term contributes 
$$
\frac{k!}{2^{k/2}(k/2)!} \sum\Sb q_1 ,\ldots ,q_{k/2} \in {\Cal P}  \\ q_i \text{ distinct} \endSb
 \ \prod_{i=1}^{k/2} \frac{h(q_i)}{q_i}\Big(1-\frac{h(q_i)}{q_i}\Big).
$$
The sum over $q$'s is bounded above by $\sigma_{\Cal P}^{k}$, and is bounded below 
by 
$$
\Big( \sum\Sb p\in {\Cal P} \\ p\ge \pi_{k/2}({\Cal P}) \endSb \frac{h(p)}{p} 
\Big(1-\frac{h(p)}{p}\Big)\Big)^{k/2} 
\ge (\sigma_{\Cal P}^2 -k/8 )^{k/2}, 
$$
where we let $\pi_{n}({\Cal P})$ denote the $n$-th smallest prime in ${\Cal P}$ 
and made use of the fact that $0\le (h(p)/p)(1-h(p)/p) \le 1/4 $.   The 
remainder of the argument is exactly the same as in Proposition 2.

Finally we need to deal with the ``error'' term contribution to (11).   To 
estimate the error terms that arise in (11), we use that
$|E(p_1\cdots p_k,m)|\leq \prod_{p_i \nmid m} h(p_i)/p_i$.  Thus the 
error term is 
$$
\leq \sum_{\ell=1}^k \sum\Sb m=q_1\dots q_\ell\geq 1 \\ q_1<q_2<\cdots <q_\ell \in \Cal P \endSb 
|r_m| \ \sum\Sb p_1,
\ldots, p_k \in \Cal P\\ m|p_1\cdots p_k\endSb \ \prod_{ p_i \nmid
m}  \frac{{h(p_i)}}{p_i}
$$
Fix $m$ and let $e_j=\#\{ i:\ p_i=q_j\}$  for each $j,\ 1\leq
j\leq \ell$. Then there are $e_0:=k-(e_1+\cdots +e_\ell)\leq
k-\ell$ primes $p_i$ which are not equal to any $q_j$, and so
their contribution to the final sum is $\leq  \mu_{\Cal P}^{e_0}$.
Therefore the final sum is
$$
 \align
 &\leq  \sum_{0\leq e_0\leq k-\ell} \binom k{e_0} \mu_{\Cal P}
^{e_0} \sum\Sb e_1+\cdots +e_\ell=k-e_0\\ \text{each}\ e_i\geq
1\endSb \frac{(k-e_0)!}{e_1!\cdots e_\ell!} \\
&\leq \sum_{0\leq e_0\leq k-1} \binom k{e_0} \mu_{\Cal P}^{e_0}
\ell^{k-e_0} \leq  (\mu_{\Cal P}+\ell)^k \ll 2\mu_{\Cal P}^k, \\
\endalign
$$
since $k^3\leq \sigma_{\Cal P}^2 \leq \mu_{\Cal P}$.   This completes the 
proof of the Proposition.

\enddemo

One way of using Proposition 3 is to take ${\Cal P}$ to 
be the set of primes below $z$ where $z$ is suitably 
small so that the error term arising from the $|r_d|$'s is 
negligible.  If the numbers $a$ in ${\Cal A}$ are not too large, 
then there cannot be too many primes larger than $z$ that divide 
$a$, and so Proposition 3 furnishes information about $\omega(a)$. 
Note that we used precisely such an argument in deducing Theorem 
1 from Proposition 2.  

In this manner, Proposition 3 may be used to prove the Erd{\H o}s-Kac theorem for 
many interesting sequences of integers.  For example, 
Halberstam [12] showed such a result for the shifted primes $p-1$, 
which the reader can now deduce from Proposition 3 and the 
Bombieri-Vinogradov theorem.  

Similarly, one can take $\Cal A=\{ f(n):\ n\leq x\}$ for
$f(t)\in \Bbb Z[t]$. In this case $h(p)$
is bounded by the degree of $f$ except at finitely many primes,
and the prime ideal theorem implies that $\mu_{\Cal P},
\sigma_{\Cal P} = m\log\log x+O(1)$ where $m$ is the number of
distinct irreducible factors of $f$.  Again this example was first considered 
by Halberstam [12].  

Alladi [1] proved an Erd{\H o}s-Kac theorem for integers without large prime 
factors.  Proposition 3 reduces this problem to obtaining information 
about multiples of $d$ in this set of "smooth numbers."  We invite 
the reader to fill in this information.  


 In place of $\omega(a)$ we may study more generally the distribution of values 
 of $g(a)$ where $g$ is an ``additive function.''   Recall that an additive 
 function satisfies $g(1)=0$, and $g(mn)=g(m)+g(n)$ whenever $m$ and $n$ are coprime.  Its values 
 are determined by the prime-power values $g(p^k)$.  If in addition $g(p^k)=g(p)$
  for all $k\ge 1$ we say that the function $g$ is 
``strongly additive."  The strongly additive functions form a particularly nice subclass of 
 additive functions and for convenience we restrict ourselves to this subclass.

\proclaim{Proposition 4} Let ${\Cal A}$ be a (multi)-set of $x$ integers, and 
let $h(d)$ and $r_d$ be as above.  Let ${\Cal P}$ be a set of primes, and 
let $g$ be a real-valued, strongly additive function with $|g(p)|\le M$ for 
all $p\in {\Cal P}$.  Let 
$$
\mu_{\Cal P}(g)= \sum_{p\in {\Cal P}} g(p)\frac{h(p)}{p}, 
\qquad \text{and} \qquad 
\sigma_{\Cal P}(g)^2 =\sum_{p\in {\Cal P}} g(p)^2 \frac{h(p)}{p}\Big(1-\frac{h(p)}{p}\Big).
$$
Then, uniformly for all even natural numbers $k\le (\sigma_{\Cal P}(g)/M)^{\frac 23}$,  
$$
\sum_{a\in {\Cal A}} \Big(\sum\Sb p|a \\ p\in {\Cal P}\endSb g(p) -\mu_{\Cal P}(g)\Big)^{k} 
= C_k x \sigma_{\Cal P}(g)^k 
\Big(1 + O\Big(\frac{k^3 M^2}{\sigma_{\Cal P}(g)^2}\Big)\Big)+O\Big(M^k \Big(\sum_{p\in {\Cal P}} \frac{h(p)}{p} \Big)^k \sum_{d\in D_k({\Cal P})} 
|r_d| \Big),
$$
while for all odd natural numbers $k\le (\sigma_{\Cal P}(g)/M)^{\frac 23}$, 
$$
\sum_{a\in {\Cal A}} \Big(\sum\Sb p|a \\ p\in {\Cal P}\endSb g(p) -\mu_{\Cal P}(g)\Big)^{k}
\ll C_k x \sigma_{\Cal P}(g)^k \frac{k^{\frac 32} M}{\sigma_{\Cal P}(g)} 
+  M^k \Big(\sum_{p\in {\Cal P}} \frac{h(p)}{p} \Big)^k \sum_{d\in D_k({\Cal P})} 
|r_d|  .
$$
\endproclaim 

 \demo{Proof} We follow closely the proofs of Propositions 2 and 3, making appropriate 
 modifications.  Let $f_r(n)$ be as 
 in the proof of Proposition 3. Then we wish to evaluate
 $$
 \sum_{a\in {\Cal A}} \Big(\sum_{p\in {\Cal P}} g(p)f_p(a)\Big)^{k} 
 = \sum_{p_1, \ldots, p_k \in {\Cal P}} g(p_1) \cdots g(p_k) 
 \sum_{a\in {\Cal A}} f_{p_1 \cdots p_k}(a).
 $$
 We may now input the results (12a,b,c) here.  Consider first the 
 error terms that arise.  Since $|g(p)|\le M$ for all $p\in {\Cal P}$ 
 this contribution is at most $M^k$ times the corresponding error 
 in Proposition 3.  To wit, the error terms are 
 $$
 \ll M^k \Big(\sum_{p\in {\Cal P}} \frac{h(p)}{p}\Big)^k \sum_{d\in D_k({\Cal P})} 
 |r_d|.
 $$
 As for the main term, note that $G(r)=0$ unless $r$ is square-full and 
 so if $q_1 < q_2 < \ldots <q_s$ are the distinct primes among the 
 $p_1$, $\ldots$, $p_k$ our main term is 
 $$
x \sum\Sb s\le k/2\endSb \sum\Sb q_1 < \ldots < q_s \\ 
 q_i \in {\Cal P} \endSb \sum\Sb \alpha_1, \ldots, \alpha_s \ge 2 \\ 
 \sum \alpha_i =k \endSb \frac{k!}{\alpha_1! \cdots \alpha_s!} \prod_{i=1}^s g(q_i)^{\alpha_i} G(q_1^{\alpha_1} \cdots q_s^{\alpha_s}). \tag{13}
 $$

  When $k$ is even there is a term with $s=k/2$ and all $\alpha_i=2$ which is 
  the leading contribution to (13).  This term contributes 
  $$
x  \frac{k!}{2^{k/2} (k/2)!} \sum\Sb q_1, \ldots, q_{k/2} \in {\Cal P}\\ q_i \text{ distinct} \endSb 
\prod_{i=1}^{k/2} g(q_i)^2 \frac{h(q_i)}{q_i} \Big(1-\frac{h(q_i)}{q_i}\Big).
$$
If we fix $q_1$, $\ldots$, $q_{k/2-1}$, then the sum over $q_{k/2}$ is 
$\sigma_{\Cal P}(g)^2 + O(M^2k)$, since $|g(p)|\le M$ for all $p\in {\Cal P}$, and 
$0\le h(p)\le p$.   Therefore the contribution of the term $s=k/2$ to (13) is 
$$
C_k x \Big(\sigma_{\Cal P}(g)^2 + O(M^2 k) \Big)^{k/2} 
= C_k x \sigma_{\Cal P}(g)^k \Big( 1+ O\Big(\frac{M^2 k^2}{\sigma_{\Cal P}(g)^2}\Big)\Big),
$$
since $kM \le \sigma_{\Cal P}(g)$.  

Now consider the terms $s<k/2$ in (13).  Since $|G(q_1^{\alpha_1} \cdots q_s^{\alpha_s})| 
\le \prod_{i=1}^{s} (h(q_i)/q_i)(1-h(q_i)/q_i)$, and $\prod_{i=1}^{s} |g(q_i)|^{\alpha_i} 
\le M^{k-2s} \prod_{i=1}^{s} |g(q_i)|^2$, we see that these terms contribute an amount 
whose magnitude is 
$$
\align
&\le x\sum_{s<k/2} \frac{k!}{s!} M^{k-2s}\Big(\sum_{q \in {\Cal P}} |g(q)|^2 \frac{h(q)}{q} 
\Big(1-\frac{h(q)}{q}\Big)\Big)^s \sum\Sb \alpha_1, \ldots, \alpha_s \ge 2\\ \sum \alpha_i =k \endSb 
\frac{1}{\alpha_1! \cdots \alpha_s!} 
\\
&\le x \sum_{s<k/2} \frac{k!}{s! 2^s} \binom{k-s}{s} M^{k-2s}\sigma_{\Cal P} (g)^{2s},
\\
\endalign
$$
using that $\binom{k-s}{s}$ equals the number of ways of writing $k=\sum \alpha_i$ with 
each $\alpha_i \ge 2$.  
The Proposition follows. 

 \enddemo 

One way to apply Proposition 4 is to take ${\Cal P}$ to be the set of 
all primes below $z$ with $|g(p)|$ being small.  If there are not too many 
values of $p$ with $|g(p)|$ being large, then we would expect that 
$g(a)$ is roughly the same as $g_{\Cal P}(a)$ for most $a$.  In such situations, 
Proposition 4 which furnishes the distribution  of $g_{\Cal P}(a)$ 
would also furnish the distribution of $g(a)$.  In this manner one can deduce 
the result of Kubilius and Shapiro [31] which is a powerful generalization of 
the Erd{\H o}s-Kac theorem for additive functions.  Indeed we can derive 
such a Kubilius-Shapiro result in the more general sieve theoretic framework 
given above, and for all additive functions rather than only for the subclass of
strongly additive functions.

There are many other interesting number theory questions 
in which an Erd\H os-Kac type theorem has been proved. 
We have collected some of these references below\footnote{Thanks
are due to Yu-Ru Liu for her help with this.} and invite the
reader to determine which of these Erd\H os-Kac type theorems 
can be deduced from the results given herein.
The reader may also be interested in the textbooks
[4,  18, 32] for a more classical discussion of some of these issues,
and to the elegant essays [2, 16].

\enddocument